# A Recursive Exponential-Gamma Mixture: A New Generalization of the Lindley Distribution

Afshin Yaghoubi, Ph.D. candidate

Faculty of Mathematics and Computer Science, Amirkabir University of Technology, Tehran, Iran;

Email: afshin.y@aut.ac.ir

Esmaile Khorram[1], Ph.D.

Faculty of Mathematics and Computer Science, Amirkabir University of Technology, Tehran, Iran;

Email: eskhor@aut.ac.ir

Omid Naghshineh Arjmand, Ph.D.

Faculty of Mathematics and Computer Science, Amirkabir University of Technology, Tehran, Iran;

Email: naghshineh@aut.ac.ir

**Abstract**

The Lindley distribution was first introduced by Lindley in 1958 for Bayesian computations. Over the past years, various generalizations of this distribution have been proposed by different authors. The generalized Lindley distributions sometimes have many parameters and, although they show good flexibility, their statistical form becomes complicated. In this study, we propose a new and simple distribution determined by the recursive relation of the Lindley distribution and the Gamma distribution with specific weights. Subsequently, some statistical properties of this distribution are examined, and finally, with several real numerical examples, its superiority over Lindley generalizations is demonstrated.

---

[1] Corresponding Author

## 1. Introduction

One of the important distributions in the field of lifetime data and reliability modeling is the Lindley distribution. Initially introduced in 1958, its statistical properties were later detailed in 2008 by [1]. They showed that this distribution exhibits better flexibility compared to the exponential distribution. Later, numerous generalizations were presented by other authors. For example, one can refer to the generalized Lindley by Zakerzadeh and Dolati [2], Nadarajah et al. [3], Shanker and Mishra [4], and Abouammoh et al. [5]. In [2], a generalization of Lindley was presented from a mixture of two gamma distributions with shape parameters $\alpha$ and $\alpha + 1$. Their proposed distribution showed high flexibility and also had a bathtub-shaped hazard rate. In [3], the cumulative distribution function (CDF) of Lindley was generalized by adding a shape parameter $\alpha$ as $F_L^\alpha(x)$. This distribution was compared with well-known gamma, log-normal, Weibull, and exponentiated exponential distributions and showed its superiority. Its hazard rate behavior also included the bathtub shape. Another distribution introduced in [4] was essentially the same as the Lindley distribution, differing only in the weights. Therefore, the authors called it the quasi-Lindley distribution. This distribution, like the Lindley, only models an increasing hazard rate. Another generalized Lindley distribution proposed by [5], resulting from a mixture of two gamma distributions with shape parameters $\alpha - 1$ and $\alpha$, also describes an increasing hazard rate behavior.

Generalizations of the Lindley distribution have recently been used in wind energy to estimate wind speed. These distributions show better fits compared to the widely used Weibull distribution, which was unrivaled in this field. For instance, the [6] introduced a 4-parameter distribution called the beta-generalized Lindley. First, its statistical properties were examined, and then it was compared with other statistical distributions for wind speed estimation, showing that their proposed distribution provided a better fit. In [7], two distributions, the generalized Lindley and power Lindley (their generalized Lindley was the same as the one introduced in [3]), were evaluated for estimating wind speed parameters. They showed that their distributions performed better than others. Another generalized Lindley distribution composed of 3 parameters was introduced by [8], demonstrating its superiority in wind speed parameter estimation compared to rival distributions. Another three-parameter generalized Lindley distribution, termed the New Alpha Power Transformed Power Lindley Distribution (NAPTPL), was proposed by [8]. The study first examined its reliability function and several related indices, and subsequently applied the model to wind speed data.

Other generalized Lindley distributions exist, such as [10-12], to which one can refer for more details.

In this article, we propose a statistical distribution. This distribution has a relatively simple mathematical form compared to generalized Lindley distributions and yet shows very good fitting power. The construction of the proposed distribution is derived from the recursive relation of the Lindley and Gamma distributions with specific weights.

The paper is organized as follows: In Section 2, the proposed distribution is introduced. In Section 3, moments, moment generating functions, and characteristic functions are examined. The Section 4, statistical indices are calculated. In Section 5, maximum likelihood estimation of the distribution parameters is presented. The distribution of the sum of $n$ i.i.d. variables and a simulation study are presented in Sections 7 and 8. In Section 9, conclusions and future work are discussed.

## 2. Description of the Proposed Distribution

In this study, we introduce a recursive generalization of the Lindley distribution obtained by combining the probability density function (PDF) of the previous stage of the Lindley distribution with scale parameter $\theta$ and the Gamma distribution with shape parameter $\alpha$ and scale parameter $\theta$. This distribution actually arises from the recursive relation:

$$f_{L_i}(x) = \frac{\theta}{\theta + 1} f_{L_{i-1}}(x; \theta) + \frac{1}{\theta + 1} f_G(x; \alpha, \theta), i = 1, 2, \ldots$$

where $f_{L_0}(x) = f_E(x; \theta) = \theta e^{-\theta x}$, $f_{L_1}(x) = \frac{\theta^2}{\theta+1}(1+x)e^{-\theta x}$, and $f_G(x; \alpha, \theta) = \frac{\theta^\alpha}{\Gamma(\alpha)} x^{\alpha-1} e^{-\theta x}$.

Applying the recursive relation for $n$ steps, the probability density function is given as follows:

$$f_X(x) = \theta e^{-\theta x} \left( p_n + (1 - p_n) \frac{(\theta x)^{\alpha-1}}{\Gamma(\alpha)} \right), x, \alpha, \theta > 0, \tag{1}$$

where $p_n = \left(\frac{\theta}{\theta+1}\right)^n$.

In fact, equation (1) is constructed from $f_X(x) = p_n f_E(x;\theta) + (1-p_n)f_G(x;\alpha,\theta)$.

When $n \to \infty$, with $(\theta > 0)$, $p_n \to 0$, in this case, (1) reduces to $f_X(x) = f_G(x;\alpha,\theta)$. For $\alpha = 1$, we obtain $f_X(x) = f_E(x;\theta)$, and for $n = 1, \alpha = 2$, $f_X(x) = f_{L_1}(x;\theta)$.

Note that in relation (1), $n$ is not a parameter but a natural number indicating the depth of the recursive generalization. By increasing $n$, the weight of the Exponential component $p_n$ decreases and the model gets closer to the Gamma distribution. In practice, we seek small $n$ values sufficient to achieve a desirable fit. Therefore, based on preliminary goodness-of-fit analyses (not reported here for brevity), we set $n = 3$ throughout the rest of this paper, as it provides a good balance between flexibility and parsimony.

Fig. 1 shows the probability density function for different parameters.

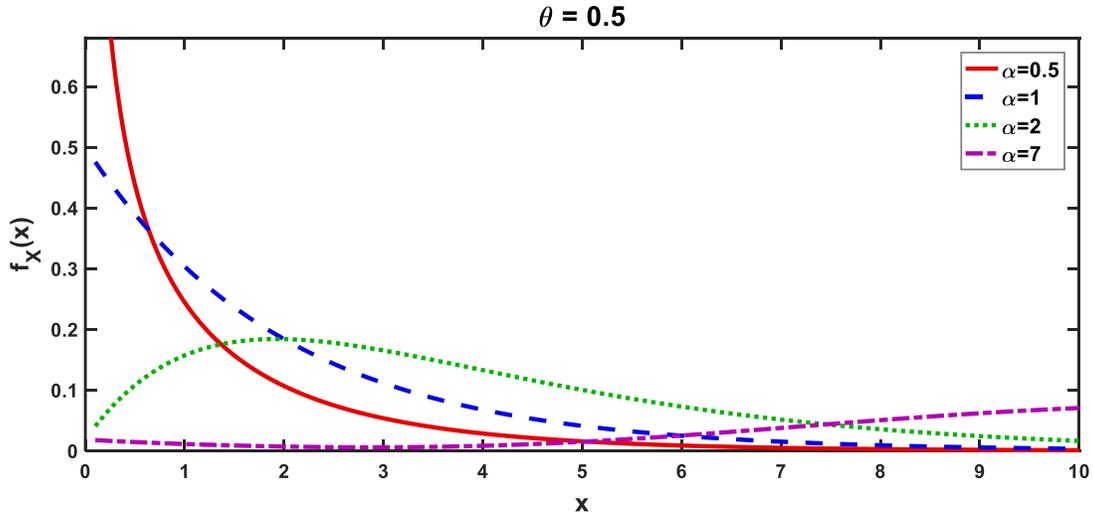

**Fig. 1.** PDF diagram

The reliability function of (1) is:

$$R_X(x) = p_n e^{-\theta x} + (1-p_n)\frac{\Gamma(\alpha,\theta x)}{\Gamma(\alpha)}, x, \alpha, \theta > 0, \qquad (2)$$

where $\Gamma(\alpha,\theta x)$ is the upper incomplete gamma function.

The CDF is also available as $F_X(x) = 1 - R_X(x)$.

The hazard function (1) is obtained from the following relation:

$$h_X(x) = \frac{f_X(x)}{R_X(x)} = \frac{\theta\left(p_n + (1-p_n)\frac{(\theta x)^{\alpha-1}}{\Gamma(\alpha)}\right)}{p_n + (1-p_n)\frac{\Gamma(\alpha,\theta x)}{\Gamma(\alpha)}e^{\theta x}}, x, \alpha, \theta > 0. \tag{3}$$

Clearly, when $\alpha = 1$, $h_X(x) = \theta$.

In special case, when the $\alpha \in N$, $\Gamma(\alpha, \theta x) = (\alpha-1)! \, e^{-\theta x} \sum_{i=0}^{\alpha-1} \frac{(\theta x)^i}{i!}$, thus $h_X(x)$ simplifies to:

$$h_X(x) = \frac{\theta\left(p_n + (1-p_n)\frac{(\theta x)^{\alpha-1}}{(\alpha-1)!}\right)}{p_n + (1-p_n)\sum_{i=0}^{\alpha-1}\frac{(\theta x)^i}{i!}}, x, \theta > 0.$$

The derivative of the hazard function; i.e., $h'_X(x)$ is:

$$h'_X(x) = \frac{\theta^{\alpha+2} x^{\alpha-2}(1-p_n)\left(-p_n\Gamma(\alpha) + (1-p_n)\left((\theta x)^{\alpha-1} - \theta x e^{\theta x}\Gamma(\alpha,\theta x)\right)\right)}{\Gamma(\alpha)\left(p_n + (1-p_n)\frac{\Gamma(\alpha,\theta x)}{\Gamma(\alpha)}e^{\theta x}\right)^2}$$

Analysis of the above derivative shows that the proposed model exhibits diverse behaviors under different parameters. It can be shown that $h_X(x)$ is decreasing for $0 < \alpha < 1$, and increasing for $\alpha = 2$ and bathtub-shaped (first decreasing then increasing) for $\alpha \geq 3$. Fig. 2 shows the hazard function plot for different parameters.

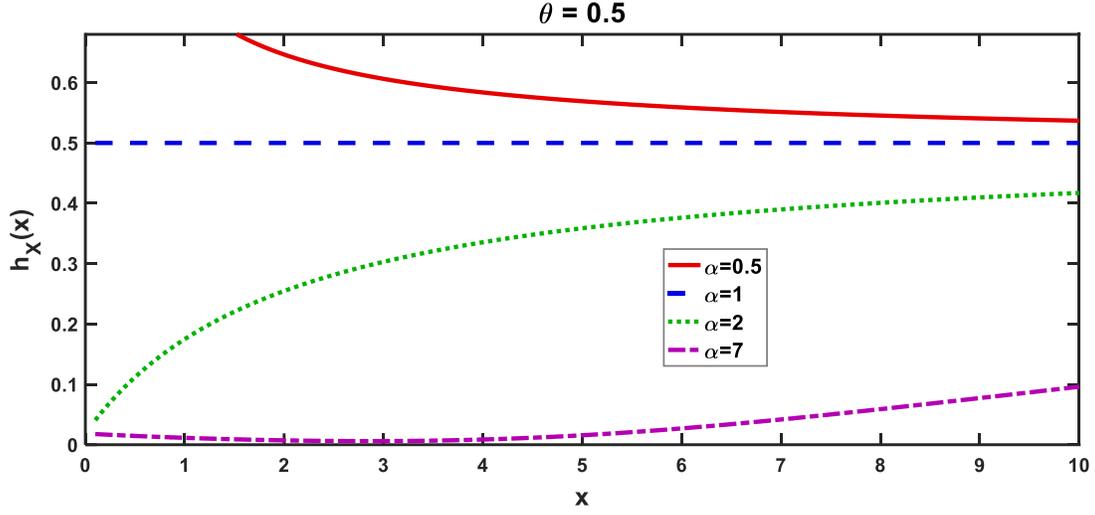

**Fig. 2.** Hazard function diagram

### 3. Moments, Moment Generating Function, and Characteristic Function

The *r*-th moment of equation (1) is:

$$\mu'_r = E(X^r) = \frac{r!}{\theta^r}\left(p_n + (1-p_n)\frac{\Gamma(\alpha+r)}{r!\,\Gamma(\alpha)}\right), r \in N,$$

and its central moment is:

$$\mu_k = E((X-\mu)^k) = \sum_{r=0}^{k}\binom{k}{r}\mu'_r(-\mu)^{k-r}, k \geq 2.$$

The moment generating function of equation (1) using the exponential series is:

$$M_X(t) = E(e^{tX}) = \sum_{r=0}^{\infty}\frac{t^r}{r!}E(X^r) = \left(1-\frac{t}{\theta}\right)^{-1}\left(p_n + (1-p_n)\left(1-\frac{t}{\theta}\right)^{-(\alpha-1)}\right), t < \theta.$$

Its characteristic function is obtained as follows:

$$\phi_X(t) = M_X(it) = \left(1-\frac{it}{\theta}\right)^{-1}\left(p_n + (1-p_n)\left(1-\frac{it}{\theta}\right)^{-(\alpha-1)}\right),$$

where $i = \sqrt{-1}$.

## 4. Statistical Indices

The first two moments of equation (1) are:

$$\mu'_1 = \frac{p_n + (1 - p_n)\alpha}{\theta},$$

$$\mu'_2 = \frac{2p_n + (1 - p_n)\alpha(\alpha + 1)}{\theta^2},$$

Now, using the obtained relations and the central moment, factors such as mean, variance, and coefficient of variation are obtained as follows:

$$E(X) = \frac{p_n + (1 - p_n)\alpha}{\theta},$$

$$Var(X) = \frac{(1 - p_n)(\alpha + p_n(1 - \alpha)^2)}{\theta^2} = \mu_2,$$

$$\gamma = \frac{\sqrt{\mu_2}}{E(X)} = \frac{\sqrt{(1 - p_n)(\alpha + p_n(1 - \alpha)^2)}}{p_n + (1 - p_n)\alpha}$$

Skewness and kurtosis are calculated from the relations $\sqrt{\beta_1} = \frac{\mu_3}{(\mu_2)^{1.5}}$ and $\beta_2 = \frac{\mu_4}{(\mu_2)^2}$.

## 5. Estimation

To estimate the parameters $\alpha$ and $\theta$, the maximum likelihood method is used. Suppose $x_1, \ldots, x_m$ are $m$ observations from distribution (1), then the likelihood function is:

$$L(\alpha, \theta) = \prod_{i=1}^{m} \theta e^{-\theta x_i} \left( p_n + (1 - p_n) \frac{(\theta x_i)^{\alpha-1}}{\Gamma(\alpha)} \right)$$

$$= \theta^m e^{-\theta \sum_{i=1}^{m} x_i} \prod_{i=1}^{m} \left( p_n + (1 - p_n) \frac{(\theta x_i)^{\alpha-1}}{\Gamma(\alpha)} \right)$$

Taking the natural logarithm of $L(\alpha, \theta)$ we have:

$$\ell(\alpha, \theta) = logL(\alpha, \theta) = mlog\theta - \theta \sum_{i=1}^{m} x_i + \sum_{i=1}^{m} log \left( p_n + (1 - p_n) \frac{(\theta x_i)^{\alpha-1}}{\Gamma(\alpha)} \right)$$

Now, by taking partial derivatives of $\ell(\alpha, \theta)$ and setting them equal to zero, we can arrive at the following system of equations:

$$\begin{cases} \dfrac{\partial \ell}{\partial \alpha} = \sum_{i=1}^{m} \left[ \dfrac{(1 - p_n) \dfrac{(\theta x_i)^{\alpha-1}}{\Gamma(\alpha)} (\log(\theta x_i) - \psi(\alpha))}{p_n + (1 - p_n) \dfrac{(\theta x_i)^{\alpha-1}}{\Gamma(\alpha)}} \right] = 0, \\ \dfrac{\partial \ell}{\partial \theta} = \dfrac{m}{\theta} - \sum_{i=1}^{m} x_i + \sum_{i=1}^{m} \left[ \dfrac{\dfrac{np_n}{\theta(\theta+1)}\left(1 - \dfrac{(\theta x_i)^{\alpha-1}}{\Gamma(\alpha)}\right) + \dfrac{(1 - p_n)(\alpha - 1)(\theta x_i)^{\alpha-1}}{\theta \Gamma(\alpha)}}{p_n + (1 - p_n) \dfrac{(\theta x_i)^{\alpha-1}}{\Gamma(\alpha)}} \right] = 0. \end{cases}$$

where $\psi(\alpha) = \dfrac{\Gamma'(\alpha)}{\Gamma(\alpha)}$ is the digamma function. The above equations can be solved using the Newton-Raphson numerical method.

### 6. Distribution of the Sum

Let $x_1, \ldots, x_m$ be independent and identically distributed (*iid*) from (1), then the distribution of the sum $S = X_1 + \cdots + X_m$ can be easily obtained using $f_X(x) = p_n f_E(x; \theta) + (1 - p_n) f_G(x; \alpha, \theta)$. We know that the moment generating function (MGF) of the Exponential distribution with scale parameter $\theta$ is equal $\left(1 - \dfrac{t}{\theta}\right)^{-1}$ and the MGF of the Gamma distribution with shape parameter $\alpha$ and scale parameter $\theta$ is $\left(1 - \dfrac{t}{\theta}\right)^{-\alpha}$. Therefore, the moment generating function of (1) is:

$$M_X(t) = p_n \left(1 - \dfrac{t}{\theta}\right)^{-1} + (1 - p_n) \left(1 - \dfrac{t}{\theta}\right)^{-\alpha}$$

Since the $X_i$ are *iid*, $M_S(t) = (M_X(t))^m$. Therefore, using the binomial expansion, $M_S(t)$ can be written as:

$$M_S(t) = \sum_{k=0}^{m} \binom{m}{k} \left(p_n \left(1 - \dfrac{t}{\theta}\right)^{-1}\right)^{m-k} \left((1 - p_n)\left(1 - \dfrac{t}{\theta}\right)^{-\alpha}\right)^{k}$$

With some simple algebraic calculations, $M_S(t)$ can be rewritten as follows:

$$M_S(t) = p_n{}^m \sum_{k=0}^{m} \binom{m}{k} \left(\frac{1}{p_n} - 1\right)^k \left(1 - \frac{t}{\theta}\right)^{-(m+k(\alpha-1))}$$

In the above equation, it is clear $\left(1 - \frac{t}{\theta}\right)^{-(m+k(\alpha-1))}$ is the MGF of a Gamma distribution with shape parameter $m + k(\alpha - 1)$ and scale parameter $\theta$. Accordingly, the density function $f_S(s)$ is:

$$f_S(s) = \sum_{k=0}^{m} \binom{m}{k} p_n{}^{m-k} (1 - p_n)^k f_G(s; m + k(\alpha - 1), \theta), s > 0.$$

## 7. Simulation Study

In this section, random data generation from (1) is performed. Since obtaining the inverse of the CDF in closed form is not feasible, the following algorithm is used to generate random numbers:

- Choose parameters $\alpha, \theta > 0$ and constant $n$,
- Calculate $p_n$,
- For each $i$, generate a random number $u$ from $U(0,1)$,
- If $u \leq p_n$, then generate a random number from the Exponential distribution; i.e., $x = -\frac{1}{\theta} \log(v), v \sim U(0,1)$,
- If $u > p_n$, then generate a random number from the Gamma distribution; i.e., $x = -\frac{1}{\theta} \sum_{j=1}^{\alpha} \log(v_j)$.

Here, we analyze the reliability of a two-component parallel system with independent and identical components, where the component failure times follow distribution (1) with parameters $(\alpha, \theta) = (3, 0.05)$, using both exact analysis and Monte Carlo simulation, and compare the results. Furthermore, in relation (1), $n = 3$ is considered.

We know that the reliability of a two-component parallel system is $R_{System}(t) = 2R_X(t) - R_X^2(t)$. Using (2), the system reliability can easily be calculated as a parametric formula in terms of $t$. Now, the system reliability is calculated using simulation. First, set the initial value $K = 0$ and generate N independent random data pairs $(X_{i1}, X_{i2})$ from the system components whose lifetimes follow distribution of (1) with $(\alpha, \theta) = (3, 0.05)$. Clearly, the lifetime for a parallel system with 2 components is $Z_i = min(X_{i1}, X_{i2})$. In order to calculate $R_{System}(t)$, $Z_i$ is compared with $t$. If $Z_i > t$, then K is increased by one. Finally, after N repetitions, $\tilde{R}_{System}(t) = \frac{K}{N}$ is calculated. Table 1 shows the results of this comparison for a time interval of 100 hours and with $N = 10,000$ repetitions.

**Table 1**

Reliability analysis

| Time ($t$) | $R_{System}(t)$ | $\tilde{R}_{System}(t)$ | Absolute Error |
|---|---|---|---|
| 0 | 1.000000 | 1.000000 | 0.000000 |
| 10 | 0.997570 | 0.997600 | 0.000030 |
| 20 | 0.958466 | 0.958500 | 0.000034 |
| 30 | 0.843260 | 0.843300 | 0.000040 |
| 40 | 0.682243 | 0.682200 | 0.000043 |
| 50 | 0.508866 | 0.508900 | 0.000034 |
| 60 | 0.356983 | 0.357000 | 0.000017 |
| 70 | 0.239365 | 0.239400 | 0.000035 |
| 80 | 0.158338 | 0.158300 | 0.000038 |
| 90 | 0.101665 | 0.101700 | 0.000035 |
| 100 | 0.063936 | 0.063900 | 0.000036 |

Based on the provided reliability analysis in Table 1, the simulation demonstrates a high degree of accuracy and internal consistency. The extremely low absolute error (consistently on the order of $10^{-5}$) between the two calculated reliability functions, $R_{System}(t)$ and $\tilde{R}_{System}(t)$, across the entire 100-hour time span serves as strong validation evidence. Furthermore, the monotonically decreasing reliability values from 1.0 to 0.064 follow a logical and expected trend for system failure over time, confirming the model's output is physically plausible and not corrupted by numerical instability or significant algorithmic errors.

## 8. Application

In this section, 4 real-world examples from different fields are presented. For accuracy and efficiency, the proposed distribution is compared with the generalized Lindley distributions of [2-5]. The probability density functions of the mentioned distributions and the proposed one are listed in Table 2:

**Table 2**

Various PDF's

| By | PDF |
|---|---|
| [2] | $f_X(x) = \dfrac{\alpha\lambda^2(1+x)\left(1 - \dfrac{1+\lambda+\lambda x}{1+\lambda}e^{-\lambda x}\right)^{\alpha-1} e^{-\lambda x}}{1+\lambda}, x, \alpha, \lambda > 0.$ |
| [3] | $f_X(x) = \dfrac{\theta^2(\theta x)^{\alpha-1}(\alpha + \gamma x)e^{-\theta x}}{(\gamma + \theta)\Gamma(\alpha + 1)}, x, \alpha, \theta, \gamma > 0.$ |
| [4] | $f_X(x) = \dfrac{\theta^\alpha x^{\alpha-2}(x + \alpha - 1)e^{-\theta x}}{(\theta + 1)\Gamma(\alpha)}, \alpha \geq 1, x, \theta > 0.$ |
| [5] | $f_X(x) = \dfrac{\theta(\alpha + \theta x)e^{-\theta x}}{\alpha + 1}, \alpha \geq -1, x, \theta > 0.$ |
| This paper | $f_X(x) = \theta e^{-\theta x}\left(p_n + (1-p_n)\dfrac{(\theta x)^{\alpha-1}}{\Gamma(\alpha)}\right), x, \alpha, \theta > 0.$ |

**Example 1.** The first real dataset relates to the waiting time (in minutes) of 100 bank customers before service. The data are given below [1]:

0.8, 0.8, 1.3, 1.5, 1.8, 1.9, 1.9, 2.1, 2.6, 2.7, 2.9, 3.1, 3.2, 3.3, 3.5, 3.6, 4.0, 4.1, 4.2, 4.2, 4.3, 4.3, 4.4, 4.4, 4.6, 4.7, 4.7, 4.8, 4.9, 4.9, 5.0, 5.3, 5.5, 5.7, 5.7, 6.1, 6.2, 6.2, 6.2, 6.3, 6.7, 6.9, 7.1, 7.1, 7.1, 7.1, 7.4, 7.6, 7.7, 8.0, 8.2, 8.6, 8.6, 8.6, 8.8, 8.8, 8.9, 8.9, 9.5, 9.6, 9.7, 9.8, 10.7, 10.9, 11.0, 11.0, 11.1, 11.2, 11.2, 11.5, 11.9, 12.4, 12.5, 12.9, 13.0, 13.1, 13.3, 13.6, 13.7, 13.9, 14.1, 15.4, 15.4, 17.3, 17.3, 18.1, 18.2, 18.4, 18.9, 19.0, 19.9, 20.6, 21.3, 21.4, 21.9, 23.0, 27.0, 31.6, 33.1, 38.5

**Example 2.** The second dataset is related to endurance tests on deep groove ball bearings, representing the number of million revolutions before failure for each of 23 bearings in life tests. The data are given below [13]:

17.88, 28.92, 33.00, 41.52, 42.12, 45.60, 48.80, 51.84, 51.96, 54.12, 55.56, 67.80, 68.44, 68.64, 68.88, 84.12, 93.12, 98.64, 105.12, 105.84, 127.92, 128.04, 173.40

**Example 3.** The third real dataset shows relief times (in minutes) for 20 patients receiving an analgesic. The data are given below [14]:

1.1, 1.4, 1.3, 1.7, 1.9, 1.8, 1.6, 2.2, 1.7, 2.7, 4.1, 1.8, 1.5, 1.2, 1.4, 3, 1.7, 2.3, 1.6, 2

**Example 4.** The fourth dataset presents the strength of aircraft window glass. The data are given below [15]:

18.83, 20.8, 21.657, 23.03, 23.23, 24.05, 24.321, 25.5, 25.52, 25.8, 26.69, 26.77, 26.78, 27.05, 27.67, 29.9, 31.11, 33.2, 33.73, 33.76, 33.89, 34.76, 35.75, 35.91, 36.98, 37.08, 37.09, 39.58, 44.045, 45.29, 45.381

For each dataset, the parameters of all models were estimated using the MLE method. The quality of fit was assessed using the Log-Liklihood (L-L), Akaike Information Criterion (AIC) and the Kolmogorov-Smirnov (K-S) test. The results are shown in Table 3:

**Table 3**

Statistical values

| Data | PDF | Estimation of Parameters | -LogL | AIC | K-S |
|---|---|---|---|---|---|
| Ex. 1 | [2] | $\hat{\alpha} = 1.156, \hat{\theta} = 0.192, \hat{\gamma} = 0.021$ | 329.45 | 662.90 | 0.0563 |
|  | [3] | $\hat{\alpha} = 1.243, \hat{\lambda} = 0.187$ | 328.92 | 663.84 | 0.0521 |
|  | [4] | $\hat{\alpha} = 1.892, \hat{\theta} = 0.204$ | 330.18 | 664.36 | 0.0614 |
|  | [5] | $\hat{\alpha} = 0.743, \hat{\theta} = 0.176$ | 331.25 | 666.50 | 0.0648 |
|  | **This paper** | $\widehat{\boldsymbol{\alpha}} = \mathbf{1.421}, \widehat{\boldsymbol{\theta}} = \mathbf{0.189}$ | **328.15** | **660.30** | **0.0487** |
| Ex. 2 | [2] | $\hat{\alpha} = 3.124, \hat{\theta} = 0.021, \hat{\gamma} = 0.005$ | 122.18 | 248.36 | 0.1834 |
|  | [3] | $\hat{\alpha} = 2.891, \hat{\lambda} = 0.031$ | 121.45 | 248.90 | 0.1762 |
|  | [4] | $\hat{\alpha} = 4.892, \hat{\theta} = 0.025$ | 123.67 | 251.34 | 0.1921 |
|  | [5] | $\hat{\alpha} = 1.892, \hat{\theta} = 0.019$ | 124.83 | 253.66 | 0.2043 |
|  | **This paper** | $\widehat{\boldsymbol{\alpha}} = \mathbf{2.743}, \widehat{\boldsymbol{\theta}} = \mathbf{0.022}$ | **120.92** | **245.84** | **0.1687** |
| Ex. 3 | [2] | $\hat{\alpha} = 2.451, \hat{\theta} = 0.892, \hat{\gamma} = 0.341$ | 24.18 | 52.36 | 0.1342 |
|  | [3] | $\hat{\alpha} = 3.124, \hat{\lambda} = 1.892$ | 23.67 | 53.34 | 0.1283 |
|  | [4] | $\hat{\alpha} = 3.892, \hat{\theta} = 1.124$ | 25.12 | 54.24 | 0.1421 |
|  | [5] | $\hat{\alpha} = 0.892, \hat{\theta} = 0.743$ | 25.89 | 55.78 | 0.1489 |
|  | **This paper** | $\widehat{\boldsymbol{\alpha}} = \mathbf{2.124}, \widehat{\boldsymbol{\theta}} = \mathbf{0.912}$ | **23.15** | **50.30** | **0.1214** |
| Ex. 4 | [2] | $\hat{\alpha} = 4.521, \hat{\theta} = 0.098, \hat{\gamma} = 0.023$ | 94.23 | 192.46 | 0.0921 |
|  | [3] | $\hat{\alpha} = 4.213, \hat{\lambda} = 0.104$ | 93.78 | 193.56 | 0.0874 |
|  | [4] | $\hat{\alpha} = 5.892, \hat{\theta} = 0.112$ | 95.12 | 194.24 | 0.0983 |

| | | | | |
|---|---|---|---|---|
| [5] | $\hat{\alpha} = 2.743, \hat{\theta} = 0.089$ | 96.34 | 196.68 | 0.1042 |
| **This paper** | $\boldsymbol{\hat{\alpha} = 3.891, \hat{\theta} = 0.101}$ | **93.25** | **190.50** | **0.0841** |

As can be seen, in Table 3 the proposed model obtained the lowest AIC value and also the lowest K-S statistic in all four datasets. This superiority indicates the flexibility and efficiency of the proposed model.

## 9. Conclusion and Future Work

In this study, a statistical distribution was introduced, constructed from the recursive relation of the Lindley and Gamma distributions with specified weights. Some of its statistical properties were stated, and a simulation study on it in the field of system reliability was performed. Furthermore, several different real-world examples were considered for the proposed distribution, and then we compared it with some generalized Lindley distributions. The results showed that the proposed distribution exhibited very good flexibility in all examples compared to other distributions. In future work, one can discuss the optimal *n*. In other words, a method for selecting the value of *n* that provides the best fit to the data can be proposed. Also, generalizing this structure to other base distributions or adding additional shape parameters could be topics for future research.

# References

[1] Ghitany, M. E., Atieh, B., & Nadarajah, S. (2008). Lindley distribution and its application. *Mathematics and computers in simulation*, *78*(4), 493-506.

[2] Zakerzadeh, H., & Dolati, A. (2009). Generalized lindley distribution.

[3] Nadarajah, S., Bakouch, H. S., & Tahmasbi, R. (2011). A generalized Lindley distribution. *Sankhya B*, *73*(2), 331-359.

[4] Shanker, R., & Mishra, A. (2013). A quasi Lindley distribution. *African Journal of Mathematics and Computer Science Research*, *6*(4), 64-71.

[5] Abouammoh, A. M., Alshangiti, A. M., & Ragab, I. E. (2015). A new generalized Lindley distribution. *Journal of Statistical computation and simulation*, *85*(18), 3662-3678.

[6] Yang, T., & Chen, D. (2025). Beta-Generalized Lindley Distribution: A Novel Probability Model for Wind Speed. *arXiv preprint arXiv:2503.09912*.

[7] Arslan, T., Acitas, S., & Senoglu, B. (2017). Generalized Lindley and Power Lindley distributions for modeling the wind speed data. *Energy Conversion and Management*, *152*, 300-311.

[8] Kantar, Y. M., Usta, I., Arik, I., & Yenilmez, I. (2018). Wind speed analysis using the extended generalized Lindley distribution. *Renewable Energy*, *118*, 1024-1030.

[9] Ahsan-ul-Haq, M., Choudhary, S. M., AL-Marshadi, A. H., & Aslam, M. (2022). A new generalization of Lindley distribution for modeling of wind speed data. Energy Reports, 8, 1-11.

[10] Algarni, A. (2021). On a new generalized lindley distribution: Properties, estimation and applications. *Plos one*, *16*(2), e0244328.

[11] Bhati, D., Sastry, D. V. S., & Qadri, P. M. (2015). A new generalized Poisson-Lindley distribution: Applications and properties. *Austrian Journal of Statistics*, *44*(4), 35-51.

[12] Ekhosuehi, N., Opone, F., & Odobaire, F. (2018). A new generalized two parameter Lindley distribution. *Journal of Data Science*, *16*(3), 549-566.

[13] Lawless, J. F. (2011). *Statistical models and methods for lifetime data*. John Wiley & Sons.

[14] Gross, A. J., & Clark, V. (1975). Survival distributions: reliability applications in the biomedical sciences. *(No Title)*.

[15]     Fuller Jr, E. R., Freiman, S. W., Quinn, J. B., Quinn, G. D., & Carter, W. C. (1994, September). Fracture mechanics approach to the design of glass aircraft windows: A case study. In *Window and dome technologies and materials IV* (Vol. 2286, pp. 419-430). SPIE.